\documentclass[12pt]{article}
\usepackage{amsfonts,amssymb,amsmath,theorem}
\newtheorem{Pa}{Paper}[section]
\newtheorem{Tm}[Pa]{{\bf Theorem}}
\newtheorem{La}[Pa]{{\bf Lemma}}
\newtheorem{Cy}[Pa]{{\bf Corollary}}
\newtheorem{Pn}[Pa]{{\bf Proposition}}
{
\theorembodyfont{\normalfont}

\newtheorem{Dn}[Pa]{{\bf Definition}}
\newtheorem{I}[Pa]{{\bf}}

}

\newcommand\abs[1]{\left|#1\right|}              %  модуль
\newcommand\norm[1]{\left\|#1\right\|}           %  норма
\newcommand\qed{\ifhmode\unskip\nobreak\fi\quad  %  квадрат в конце
   \ifmmode\square\else\hbox{$\square$}\fi}      %   доказательств

\let\cal=\mathcal

\newcommand{\Nc}{\mathcal{N}}

\newcommand{\al}{\alpha}

\newcommand{\End}{\operatorname{End}}
\newcommand{\A}{\mathcal A}

\numberwithin{equation}{section}

\begin{document}
\begin{center}
{\bfseries \large Crossed product of a  $C^*$-algebra by
a semigroup of endomorphisms generated by partial isometries}

\bigskip
B.K. Kwa\'sniewski, A.V. Lebedev

\bigskip
 University of Bialystok;\\
 University of Bialystok / Belarus State University
\end{center}

\vspace{5mm}
\quad\parbox{14.5cm}{\small \hspace{0.5cm}
The paper presents a
construction of the crossed product of a  $C^*$-algebra by a semigroup of endomorphisms generated by partial isometries.

\medbreak

\medbreak
 }

\vspace{5mm}
\tableofcontents

\section{Introduction}\label{intro}

Given a $C^*$-algebra $\A$ and an endomorphism $\alpha$ there is a
number of ways to construct a new $C^*$-algebra (an extension of
$\A$) called the crossed product. Among the successful
constructions of this sort one would mention, for example,  the
constructions developed by J. Cuntz and  W. Krieger
\cite{cuntz,CK}, W.L. Paschke \cite{Paschke},
 P.J. Stacey \cite{Stacey}, G.J. Murphy \cite{Murphy}, R. Exel
 \cite{exel}, and
B.K. Kwasniewski \cite{kwa}. Recently A.B. Antonevich, V.I Bakhtin and A.V. Lebedev
have developed in \cite{Bakht-Leb,Ant-Bakht-Leb} the construction of the crossed product that unifies
all the previous structures in the situation when the endomorphism is generated by a partial isometry.
In addition the crossed product elaborated possesses 'almost all' the fine properties of
the crossed products associated with automorphisms: being started with the semigroup ${\mathbb Z}^+$ generated by $\alpha$ it
can be represented by means of Fourier type serieses over the 'group' ${\mathbb Z}$ (associated with $\alpha$ and the corresponding transfer operator), it has a reasonable regular representation and its faithful representations are described by Isomorphism Theorem (Theorem 3.6, \cite{Ant-Bakht-Leb}) which is a quite appropriate analogy to the corresponding results for the automorphisms case (cf. Remark 3.8 \cite{Ant-Bakht-Leb}). This situation calls natural remenicences of say partial automorphisms situation
where R.Exel \cite{Exel1} invented the crossed product by a {\em single} partial automorphism and then K.  McClanahan
\cite{McClanahan} developed this construction further up to partial actions of an {\em arbitrary group}. His results
along with those obtained by A.V. Lebedev \cite{top-free} show that the crossed product associated with partial action
 in fact behaves  ideally  --- almost like in the situations of the action by automorphisms. This reminiscence is one of the main impulses
 that caused the appearance of the present article ---  once the construction of the crossed product generated by a {\em single}
  endomorphism appeared in \cite{Ant-Bakht-Leb}  it is a natural desire to develop it from the situation of a {\em single} endomorphism
  to a {\em semigroup} of endomorphisms. Recalling also that in \cite{Ant-Bakht-Leb} there sprang out a natural passage from the semigroup
   ${\mathbb Z}^+$  to the 'group' ${\mathbb Z}$ one feels a spontaneous wish to associate with a given semigroup of endomorphisms a certain
   'group' generated by endomorphisms and the corresponding transfer operators. The fulfillment of this wish is the main theme of the article.

It should by noted at once that the situation with endomorphisms is {\em not} quite the same as the situation with partial automorphisms
since the endomorphism situation is 'heavily' irreversible. The naturally arising here in the crossed product construction
notion of a {\em finely representable pair} (see 2.4, \cite{Ant-Bakht-Leb}) and the notion of a {\em finely representable
 $C^*$-dynamical system} (see \ref{fine} of the present article) 'feels' what is 'left' and what is 'right' (see (\ref{c11}) and
 Proposition \ref{universality}). Involving the physical associations one can also say that it feels the 'past' and the  'future'.
 In other words it feels
the {\em order}. Therefore it seems that the natural development of the group $\mathbb Z$ case here is the case of a {\em totally ordered}
group. And the material of the paper shows that this case can be worked out perfectly.

We  also note that a special case of the crossed product elaborated in the present article was considered by
J. Lindiarni and I. Raeburn in \cite{Lin-Rae}.

The paper is organized as follows.

In Section \ref{the_system} we introduce   the necessary notion of an action $\al$ of the positive cone
$\Gamma^+$ of a totally ordered abelian group $\Gamma$ by endomorphisms of a $C^*$-algebra $\A$  thus defining  a
$C^*$-dynamical system  $(\A,\Gamma^+,\al)$. By using the ideas and methods that originated in essence in \cite{exel,Bakht-Leb}
we  deduce  the existence of a complete transfer action for  $(\A,\Gamma^+,\al)$  (Theorem \ref{complete}).
In Section \ref{endom} we discuss the principal  constructive element of the crossed product defined further --- the finely representable
$C^*$-dynamical systems. Here the results form a natural development of the corresponding results of \cite{Bakht-Leb}.
Section \ref{l1}
is devoted to the explicit presentation of  the 'base' of the crossed product --- the Banach $^*$-algebra $l_1(\Gamma,\al,\A)$.
The crossed product $\A \rtimes_\al \Gamma$ itself is then  defined in Section \ref{Cr} and the main structural result here is
Theorem \ref{crossedstructure}. Finally in Section \ref{Faith}  we give a criterion for a representation of  $\A \rtimes_\al \Gamma$
to be faithful and present the regular representation of $\A \rtimes_\al \Gamma$.

\section{The dynamical system $(\A, \Gamma^+, \alpha)$ and
\\
transfer actions}\label{the_system}
The results of this section form in essence a natural development of the corresponding results from \cite{Bakht-Leb,exel}
onto the totally ordered abelian group situation.

Let  $\cal A$ be a $C^*$-algebra with an identity $1$ and let
$\Gamma^+$ be the positive cone of a totally ordered abelian group
$\Gamma$ with an identity $0$:
$$
\Gamma^+=\{x\in \Gamma: 0\leq x\}, \qquad \Gamma=\Gamma^+ - \Gamma^+.
$$
We  fix a  semigroup homomorphism $\al:\Gamma^+\rightarrow \End(\A)$, that is
$$
  \al_0={\rm Id}, \quad \qquad  \al_x\circ \al_y =\al_{x+y}, \quad\qquad \al_x,\,\,\al_y\in \End(\A),\,\,\, x,y\in\Gamma^+.
$$
Depending on inclinations one may say that  we fix  an action  $\alpha$  of $\Gamma^+$ by $^*$-endomorphisms of $\A$, or  a
\emph{$C^*$-dynamical system} $(\A,\Gamma^+,\al)$.
 In the sequel  we shall often make use of the simple fact that $\{\al_x(1)\}_{x\in \Gamma^{+}}$ is a nonincreasing family of projections.
 Indeed,  $\al_x(1)$ are self-adjoint idempotents and if $x\leq y$, that is $y-x\in\Gamma^+$, then   $$
\al_x(1)\al_y(1) =\al_{x}(1)\al_{x}(\al_{y-x}(1))
=\al_{x}( 1 \al_{y-x}(1))= \al_{x+y-x}(1)=\al_y(1)
,$$
that is $x\leq y \Rightarrow \al_x(1)\geq \al_y(1)$.

Throughout the article as a rule we shall denote by $x,y$ the elements of $\Gamma^+$ and by $g$ an element of $\Gamma$.

 \begin{I} Let  $L$ be an action  of $\Gamma^+$ by continuous, linear,
positive maps
$L_x:\A\to \A$, that is
$$
 L_{x+y}=L_y\circ L_x, \qquad x,y\in\Gamma^+.
$$
We shall say that $L$ is a \textit{transfer action} or an
\textit{action by transfer operators} for $(\A,\Gamma^+,\al)$, if
the following identity is satisfied
\begin{equation}\label{b,,2}
L_x(\alpha_x(a)b)=aL_x(b)
\end{equation}
for all $a,b\in \A, x\in \Gamma^+$. Note that then  $L_x(b\alpha_x(a))=L_x(b)a$ as well.
\\
These  relations have the following consequences (see \cite{exel}):  $L_x(\A)$ is a two-sided ideal, $L_x(1)$  is  a positive central
element in $\A$ (and therefore   $L_x(1)\A$ is a two-sided ideal), and the following formula holds
\begin{equation}
\label{1}
L_x(\,\cdot\,) = L_x(\al_x(1)\,\cdot\,) =
L_x(\,\cdot\,\al_x(1)).
\end{equation}
\end{I}
\begin{I}
\label{nondegeneracy}
 The  transfer action  $L$ will be called {\em non-degenerate} if for each $x\in\Gamma^+$ one of  the equivalent
  conditions holds (see \cite{exel}, Proposition 2.3):

\smallskip
\quad\ \llap{$(i)$}\ \ the composition $E_x = \al_x \circ L_x$
is a conditional expectation onto $\al_x(\A)$,

\smallskip
\quad\ \llap{$(ii)$}\ \ $\al_x\circ L_x\circ\al_x =\al_x$,

\smallskip
\quad\ \llap{$(iii)$}\ \ $\al_x (L_x(1)) = \al_x (1)$.
\bigskip

\noindent In particular, as $\al_0={\rm Id}$ (ii) implies  $L_0={\rm Id}$.
\bigskip

\noindent The non-degeneracy of the transfer action  implies (see \cite{Bakht-Leb}, Propositions 2.5, 2.6)
 that: $L_x(\A)=L_x(1)\A$,  the element $L_x(1)$ is a  central orthogonal projection in $\cal A$,  and  $L_x \!: \al_x( \A) \to L_x
(\A)$ is a $^*$-isomorphism where the inverse one is  $\al_x \!:
L_x(\A) \to \al_x(\A)$.  Moreover we note that $\{L_x(1)\}_{x\in
\Gamma^{+}}$ is a nonincreasing family of projections. Indeed, if
$x\leq y$, that is $y-x\in\Gamma^+$, then
 $$
L_x(1)L_y(1) =L_{x}(1)L_{x}\big(L_{y-x}(1)\big)
=L_{x}\big(\al_x(L_{x}(1))L_{y-x}(1)\big)= L_{x}\big(\al_x(1)L_{y-x}(1)\big)
$$
$$
= L_{x}(L_{y-x}(1))=L_y(1)
,$$
that is $x\leq y \Rightarrow L_x(1)\geq L_y(1)$.
\end{I}
\begin{I}\label{complete transfer action}
The transfer action $L$ will be called {\em complete}, if
\begin{equation}\label{completeness property}
\al_x(L_x(a))=\al_x(1)a\al_x(1),\quad \qquad x\in \Gamma^+,\quad a\in \A.
\end{equation}
The completeness  of the transfer action  implies
 that for each $x\in \Gamma^+$ we have $
 \al_x\circ L_x\circ \al_x=\al_x$ and  recalling \eqref{1} we also obtain  $L_x\circ \al_x\circ L_x= L_x
 $. In particular, a complete transfer action is non-degenerate, in addition if $L$ is complete then for
 each $x\in \Gamma^+$, we have $\al_x(\A)=\al_x(1)\A\al_x(1)$  that is $\al_x(\A)$ is a hereditary subalgebra of
  $\A$ (see \cite{exel}, Proposition 4.1).
\end{I}
The existence and uniqueness of a transfer action for an arbitrary
dynamical system $(\A,\Gamma^+,\al)$ is quite a problematic
matter. However, thanks to \cite{Bakht-Leb} one can write down
conditions on $(\A,\Gamma^+,\al)$ under which  a unique complete
transfer action for $(\A,\Gamma^+,\al)$ does exist. The next
result is a development of Theorem 2.8, \cite{Bakht-Leb} onto the
situation under investigation.
\begin{Tm}\label{complete}
Let $(\A,\Gamma^+,\al)$  be  a dynamical system. The following are equivalent:
\begin{itemize}
\item[$1)$]   there exists a complete transfer action  $L$   for $(\A,\Gamma^+,\al)$,
\item[$2)$] $(i)$ there exists a non-degenerate transfer action  $L$   for $(\A,\Gamma^+,\al)$,\\[6pt]
(ii)   $\al_x(\A)$ is   hereditary   subalgebra of $\A$ for each $x\in\Gamma^+$,
\item[$3)$] $(i)$ there exists a  family $\{P_x\}_{x\in\Gamma^+}$ of  central orthogonal projections in $\A$ such that
\begin{itemize}
\item[$a)$] $\al_x(P_{x+y}) = \al_x (1)P_y$, for all $x,y\in \Gamma^+$,
\item[$b)$] the mappings $\al_x \!: P_x \A \to
\al_x (\A)$ are $^*$-isomorphisms, and
\end{itemize}
$(ii)$ $\al_x(\A)=\al_x(1)\A\al_x(1)$ for each $x\in\Gamma^+$.
\end{itemize}
Moreover the objects in $1)$ -- $3)$ are defined in a unique way
$($i.\,e.\ the transfer action $L$ in $1)$ and $2)$ is
unique and the family  of projections $\{P_x\}_{x\in\Gamma^+}$ in $3)$ is unique as well$)$ and
\begin{equation}
\label{P} P_x=L_x(1),\qquad x\in\Gamma^+,
\end{equation}
and
\begin{equation}\label{d*-}
L_x(a) =\al_x^{-1}(\al_x(1)a\al_x(1)), \ \ \  a\in \A, \ \
\end{equation}
where $\al_x^{-1}:\al_x(\A)\to P_x\cal A$ is the inverse
mapping to\ \ $\al_x:P_x\A \to \al(\A)$,  $x\in \Gamma^+$.
\end{Tm}
{\bfseries Proof.} $1)\Rightarrow 2)$. Follows from the definition of a complete transfer action, see \ref{complete transfer action}.
\\
$2)\Rightarrow 3)$. It is known that $2)$ $(ii)$ and $3)$ $(ii)$ are equivalent, cf. \cite{exel}, Proposition 4.1.\\
 Set $P_x=L_x(1)$, $x\in \Gamma^+$. By \ref{nondegeneracy}, $\{P_x\}_{x\in\Gamma^+}$ is a family of  central orthogonal
 projections and  the mappings  $\al_x \!: P_x\A \to \al_x(\A)$ are  $^*$-isomorphisms. So $3)$ $(i)$ $b)$ is true. Recalling
 (\ref{1}) and using the fact that   $L_x \!: \al_x(A) \to P_x\A$ is the inverse to $\al_x \!: P_x\A \to \al_x(\A)$ we obtain
  $$
 \al_x(P_{x+y})=\al_x(L_{x+y}(1))=\al_x(L_{x}(L_y(1)))=\al_x\big(L_{x}(\al_x(1)P_y\al_x(1))\big)
 $$
 $$
 =\al_x(1)P_y\al_x(1)=\al_x(1)P_y,
 $$
 which proves  $3)$ $(i)$ $a)$.\\
$3) \Rightarrow 1)$. Fix $x\in \Gamma^+$. Let $\al_x^{-1}\!:\al_x(\A)\to
P_x\cal A$ be the inverse mapping to \ $\al_x\!:P_x\cal A\to
\al_x(\A)$. Define the operator  $L_x$ by the formula
$L_x(a) =\al_x^{-1}(\al_x(1)a\al_x(1))$. Clearly
$L_x$ is linear and  positive, and (\ref{completeness property}) is fulfilled. Note that \begin{equation*}
\al_x\big(L_x(\al_x(a)b)\big) =\al_x(1)\al_x(a)b
\al_x(1) =\al_x(a)\al_x(1)b\al_x(1) =\al_x(a L_x(b)),
\end{equation*}
and as the elements  $L_x(\al_x(a)b)$
and
$aL_x(b)$ belong to the ideal  $P_x\A$ where the endomorphism $\al_x$
is injective,  they coincide. Therefore \eqref{b,,2} holds, and the only thing left to  prove is that
 $L$ is an action of the semigroup $\Gamma^+$.
\\
For that purpose let us observe that the family $\{P_x\}_{x\in\Gamma^+}$ is nonincreasing. Indeed,
if $x,y\in \Gamma^+$ are such that $x\leq y$, then   $\A$ can be written as the direct sum of ideals
in two ways $\A=\ker\al_x \oplus (P_x\A)=\ker\al_y\oplus (P_y\A)$, and $\ker \al_x\subset \ker\al_y$, whence $P_{y}\A\subset P_x\A$
 and hence $P_y\leq P_x$.
\\
Using $3)$ $(i)$ $a)$,  we have  $
\al_y(P_{x+y}\A)=\al_y(P_{x+y})\al_{y}(\A)= P_x\al_y(\A)$   and as $P_{x+y}\A\subset P_y\A$ we obtain  that
$\al_y:P_{x+y}\A\rightarrow P_x\al_y(\A)$ is a $^*$-isomorphism and the inverse is given by $L_y$.  Thus    we have
$$
L_y(L_x(\A))=L_y(P_x\A)=L_y(\al_y(1)P_x\A\al_y(1))=L_y(P_x\al_y(\A))=P_{x+y}\A.
$$
Hence $L_y(L_x(a)$ and $L_{x+y}(a)$ belong to the ideal  $P_{x+y}\A$ where the endomorphism $\al_{x+y}$
is injective, and  as
$$
\al_{x+y}(L_y(L_x(a))=\al_{x}\big(\al_y(L_y(L_x(a))\big)=\al_{x}\big(\al_y(1)L_x(a)\al_y(1)\big)=\al_{x+y}(1)\al_x(L_x(a))\al_{x+y}(1)
$$
$$
=\al_{x+y}(1)a\al_{x+y}(1)=\al_{x+y}(L_{x+y}(a))
$$
we have  $L_{x+y}=L_y\circ L_x$.

The uniqueness of the objects in 1) - 3) can be established
absolutely in the same way as it is done in the proof of Theorem
2.8 in \cite{Bakht-Leb}. \qed

\section{Finely representable systems}
\label{endom}
Here we discuss one of the main notions of the article --- the finely representable systems. They  play a
principal role as in the construction of the  crossed product in Section \ref{Cr}
so also in the construction of its  regular representation in Section  \ref{Faith}.

\begin{I}
\label{fine}
Let $\cal A$  be a $C^*$-algebra  containing an identity $1$ and
$\al\!: \Gamma^+\to End(\cal A)$ be a semigroup homomorphism.
We say that the triple $(\A,\Gamma^+,\al)$ is
{\em finely representable\/} if  if
there exists a triple $(H,\pi ,U)$ consisting of a Hilbert space
$H$, faithful non-degenerate representation $\pi\!:{\cal A}\to
L(H)$ and a semigroup homomorphism $U\!:\Gamma^+ \to L(H)$   such that
 for every $a\in\A$, $x\in\Gamma^+$, the
following conditions are satisfied
\begin{gather}\label{b,,4}
\pi(\al_x(a)) =U_x\pi(a)U^*_x,\qquad U^*_x\pi(a)U_x \in \pi(\A)\\
\intertext{and}
\label{c11}
U_x\pi (a) = \pi(\al_x(a))U_x, \qquad a \in \cal A.
\end{gather}
In this case we also say that $\A$ is a {\em coefficient algebra
associated with $\al$.}

{\bf Remark.} The notion of a coefficient algebra associated with a {\em single} endomorphism was introduced
by A.V. Lebedev and A. Odzijewicz in
\cite{Leb-Odz} and proved to be of principle importance in the investigation of extensions of $C^*$-algebras by partial
isometries. We shall also use certain ideas and methods from \cite{Leb-Odz} while uncovering the internal structure of
the crossed product  in Sections \ref{Cr} and \ref{Faith}.
\medskip

It is an easy exercise that in the above definition one can replace  condition \eqref{c11} by
the condition
\begin{equation}
\label{c1*} U_x^*U_x\in Z(\pi ({\cal A})), \ \ \ x\in\Gamma^+
\end{equation}
where $Z(\pi ({\cal A}))$ is the centrum of $\pi ({\cal A})$,  or the condition
\begin{equation} \label{c9*}
 U_x^*U_x\pi (a)U_x^* \pi (b)U_x = \pi (a)U_x^*\pi (b)U_x, \ \ a,b \in \cal A, \ \ \ x\in \Gamma^+ .
\end{equation}

In particular it is clear that the finely representable $C^*$-dynamical system   can also be defined
 as a triple  $(\A,\Gamma^+,\al)$
  such that there exists a
triple $(H,\pi ,U)$ where $\pi\!:{\cal A}\to L(H)$ is a faithful
non-degenerate representation, $U_x\in L(H), \ x\in \Gamma^+$ and the mappings $U_x\,\cdot\,
U_x^*, \ \ x\in \Gamma^+$ coincide with the endomorphisms $\alpha_x$  on $\pi (\cal A)$ while
the mappings $U_x^*\,\cdot\,U_x $ are  transfer operators for $\alpha_x$.
\end{I}

The next theorem presents a criterion of the fine representability. The construction given
in its proof is in fact a development of the corresponding construction from Theorem 3.1, \cite{Bakht-Leb} onto the situation under consideration.

\begin{Tm}\label{b..1}
$(\A,\Gamma^+,\al)$ is finely representable  iff there exists a complete transfer action $L$
for $({\cal A}, \Gamma^+, \al)$ (that is either of the equivalent
conditions of Theorem \ref{complete} hold).
\end{Tm}

{\bfseries Proof.} \, {\em Necessity}. \, If conditions (\ref{b,,4}) and
(\ref{c11})  are satisfied then (identifying $\cal A$ with
$\pi (\cal A)$) one can set
\[
L_x (\cdot ) = U_x^* (\cdot )U_x, \qquad x\in \Gamma^+
\]
and it is easy to verify that $L$ is a complete transfer
action.

\smallbreak
{\em Sufficiency.} \, Let $L$ be a complete transfer action.
We shall construct the desired Hilbert space $H$ by means of the
elements of the initial algebra $\cal A$ in the following way. Let
$\langle\,\cdot\,,\, \cdot\,\rangle$ be a certain non-negative inner
product on $\cal A$ (differing from a common inner product only in
such a way that for certain non-zero elements $v\in\cal A$ the
expression $\langle v,v\rangle$ may be equal to zero). For example
this inner product may have the form $\langle v,u\rangle =f(u^*v)$
where $f$ is some positive linear functional on $\cal A$. If one
factorizes $\cal A$ by all the elements $v$ such that $\langle
v,v\rangle =0$ then he obtains a linear space with a strictly positive
inner product. We shall call the completion of this space with respect
to the norm $\norm v =\sqrt{\langle v,v \rangle}$ the {\em Hilbert
space generated by the inner product\/} $\langle\,\cdot\,,\,
\cdot\,\rangle$.

Let $F$ be the set of all positive linear functionals on $\A$. The
space $H$ will be the  direct sum $\bigoplus_{f\in
F}H^f$ of some Hilbert spaces $H^f$. Every $H^f$ will in turn be the
 direct sum of Hilbert spaces $\bigoplus_{g\in \Gamma} H^f_g$. The spaces $H^f_g$ are generated by non-negative inner
products $\langle\,\cdot\,,\,\cdot\,\rangle_g$ on the initial algebra
$\cal A$ that are given by the following formulae
\begin{align}
\label{b,,7}
\langle v,u\rangle_0 &=f(u^*v);\\[3pt]\label{b,,8}
\langle v,u\rangle_x &=f\bigl(L_x(u^*v)\bigr),
\qquad x \in \Gamma^+;\\[3pt]
\label{b,,9}
\langle v,u\rangle_{-x}
&=f\bigl(u^*\al_x(1)v\bigr), \qquad x\in \Gamma^+.
\end{align}
The properties of these inner products are described in  the next

\begin{La}
\label{b..2}
For any $v,u\in\cal A$, any $x \in \Gamma^+$ and any $g\in \Gamma$ the following equalities are true
\begin{align}
\label{b,,10} \langle \al_x(v),u\rangle_{g} &=\langle
v,L_x(u)\rangle_{g-x}, \qquad x\leq g;\\[3pt]
\label{b,,10.5} \langle \al_x(1)\al_g(v),u\rangle_{g} &=\langle
v,L_g(u)\rangle_{g-x}, \qquad  0\leq g\leq x;\\[3pt]
\label{b,,11} \langle \al_{x-g}(1)v,u\rangle_{g} &=\langle
v,\al_{x-g}(1)u \rangle_{g-x}, \,\,\,\,  g\leq 0.
\end{align}
\end{La}

{\bfseries Proof.}
Let $x\leq g$. The proof of  (\ref{b,,10}) reduces to the verification of the equalities
\[
L_g(u^*\al_x(v))=L_{g-x}\big(L_x(u^*\al_x(v)\big)=L_{g-x}\big(L_{x}(u)^* v\big)
\]
which follow from  the definition of the transfer action (recall $x\leq g$ iff $g-x\in \Gamma^+$).\\
Let $0\leq g \leq x$. Formula
(\ref{b,,10.5}) follows from
$$
L_g\big(u^*\al_x(1)\al_g(v)\big)=L_g\big(u^*\al_g(\al_{x-g}(1))\al_g(v)\big)=L_g\big(u^*\al_g(\al_{x-g}(1) v)\big)
 =L_y(u)^* \al_{x-y}(1) v
$$
Let $g\leq 0$. The proof of
(\ref{b,,11}) reduces to the verification of the equality
\[
u^*\al_{-g}(1)\al_{x-g}(1)v = u^*\al_{x-g}(1)v
\]
which follows from $\al_{-g}(1)\geq \al_{x-g}(1)$.
\qed

\smallbreak
Now let us define the  semigroup homomorphism $U\!:\Gamma^+ \to L(H)$. For an arbitrary fixed $x\in \Gamma^+$
the operators  $U_x$ and $U^*_x$ will  leave invariant all the subspaces
$H^f\subset H$. The action of  these operators  on every $H^f=\bigoplus_{g\in \Gamma} H^f_g$ is similar
 and its scheme is presented in the the next
diagrams.
\\
\setlength{\unitlength}{0.95mm}
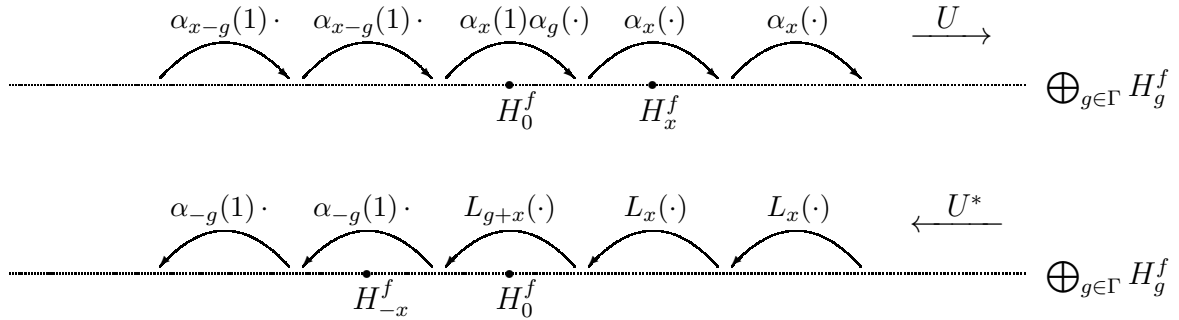
\centerline{
\begin{picture}(80,27)(-37,-9)
\qbezier[300](-80,0)(-15,0)(62,0)\put(65,-1.5){$\bigoplus_{g\in \Gamma} H^f_g $}
\put(46,8){$\underrightarrow{\,\,\,\,\,U\,\,\,\,\,\,}$}
\put(10,0){\circle*{1}}\put(8,-5){$H_{x}^f$}
\put(-10,0){\circle*{1}}\put(-12,-5){$H_0^f$}
\small
\qbezier(39,)(30,11)(21,1) \put(38,2){\vector(1,-1){1}}\put(26,8){$\al_x(\cdot)$}
\qbezier(19,)(10,11)(1,1) \put(18,2){\vector(1,-1){1}}\put(6,8){$\al_x(\cdot)$}
\qbezier(-1,1)(-10,11)(-19,1) \put(-2,2){\vector(1,-1){1}}\put(-17,8){$\al_{x}(1)\al_{g}(\cdot)$}
\qbezier(-21,1)(-30,11)(-39,1) \put(-22,2){\vector(1,-1){1}}\put(-37.5,8){$\al_{x-g}(1)\,\cdot$}
\qbezier(-41,1)(-50,11)(-59,1) \put(-42,2){\vector(1,-1){1}}\put(-57.5,8){$\al_{x-g}(1)\,\cdot$}
\end{picture}}
\\
\setlength{\unitlength}{0.95mm}
\centerline{
\begin{picture}(80,27)(-37,-10)
\qbezier[300](-80,0)(-15,0)(62,0)\put(65,-1.5){$\bigoplus_{g\in \Gamma} H^f_g $}
\put(46,8){$\underleftarrow{\,\,\,\,\,\,\,U^*\,\,\,\,}$}
\put(-30,0){\circle*{1}}\put(-32,-5){$H_{-x}^f$}
\put(-10,0){\circle*{1}}\put(-12,-5){$H_0^f$}
\small
\qbezier(39,)(30,11)(21,1) \put(22,2){\vector(-1,-1){1}}\put(26,8){$L_x(\cdot)$}
\qbezier(19,)(10,11)(1,1) \put(2,2){\vector(-1,-1){1}}\put(6,8){$L_x(\cdot)$}
\qbezier(-1,1)(-10,11)(-19,1) \put(-18,2){\vector(-1,-1){1}}\put(-16.5,8){$L_{g+x}(\cdot)$}
\qbezier(-21,1)(-30,11)(-39,1) \put(-38,2){\vector(-1,-1){1}}\put(-37.5,8){$\al_{-g}(1)\,\cdot$}
\qbezier(-41,1)(-50,11)(-59,1) \put(-58,2){\vector(-1,-1){1}}\put(-57.5,8){$\al_{-g}(1)\,\cdot$}
\end{picture}}
\noindent
Formally this action is defined in the
following way. For any \emph{finite} sum
$$
h =\bigoplus_g h_g \in H^f,\qquad  h_g\in H^f_g,
$$
we set
$$
U_x h =\bigoplus_g(U_xh)_g
\quad  \textrm{ and } \quad
U_x^* h =\bigoplus_g(U_x^*h)_g
$$
where
$$
(U_x h)_g=
\begin{cases}
\al_x(h_{g-x}),&\ \ {\rm if}\ \ x\leq g,\\
\al_x(1) \al_g(h_{g-x}),&\ \ {\rm if}\ \ 0 \leq g\leq x,
\\
\al_{x-g}(1)h_{g-x},& \ \ {\rm if}\ \ g \leq 0,
\end{cases}
$$
$$
(U^*_xh)_g =
\begin{cases}
L_x(h_{g+x}), & \ \ {\rm if}\ \ 0\leq g,\\
L_{g+x}(h_{g+x}), & \ \ {\rm if}\ \ -x\leq g\leq 0,\\
\al_{-g}(1)h_{g+x},&\ \  {\rm if}\ \ g <0.
\end{cases}
\label{b,,13}
$$
Lemma  \ref{b..2} guarantees that the operators  $U_x$ and $U^*_x$ are
well defined (i.\,e.\ they preserve factorization and completion by
means of which the spaces $H^f_g$ were built from the algebra~$\cal A$)
and $U_x$ and $U^*_x$ are mutually adjoint.
\par
Let us show that $U$ is  a semigroup homomorphism. Take any  $x,y\in \Gamma^+$ and  $h_g \in H^f_g$.\\
For $x+y\leq g$ we have $y\leq g$ and $x\leq g-y$, whence
$$
\big(U_y(U_x h)\big)_g=\al_y\big((U_xh)_{g-y}\big)=\al_y(\al_x(h_{g-y-x})=\al_{x+y}(h_{g-(x+y)})=(U_{x+y}h)_g.
$$
For $0\leq g\leq x+y$ the two cases are possible:
If $y\leq g$ then we have  $0\leq g-y\leq x$, whence
$$
\big(U_y(U_x h)\big)_g=\al_y\big((U_xh)_{g-y}\big)=\al_y\big(\al_x(1) \al_{g-y}(h_{g-y-x})\big)=\al_{x+y}(1)\al_g(h_{g-(x+y)})=(U_{x+y}h)_g.
$$
If $g\leq y$ then we have  $ g-y\leq 0$, and thus
$$
\big(U_y(U_x h)\big)_g=\al_y(1)\al_g\big((U_xh)_{g-y}\big)=\al_y(1)\al_{x+y}(1)\al_g(h_{g-(x+y)})=(U_{x+y}h)_g
$$
where in the final equality we used the fact that $\al_{x+y}(1)\leq \al_y(1)$.
\\
For $g\leq 0$ we have
$$
\big(U_y(U_x h)\big)_g=\al_{y-g}(1)(U_xh)_{g-y}=\al_{y-g}(1)\al_{x+y-g}(1)h_{g-(x+y)}=(U_{x+y}h)_g.
$$
where in the final equality we used the inequality $\al_{y-g}(1)\geq \al_{x+y-g}(1)$. Thus we have  proved
that $U\!:\Gamma^+ \to L(H)$ is a semigroup homomorphism.
\par
Now let us define the representation  $\pi\!:{\cal A}\to L(H)$.
For any $a\in\cal A$ the operator  $\pi(a)\!:H\to H$ will leave
invariant all the subspaces $H^f\subset H$ and also all the
subspaces  $H^f_g \subset H^f$. If  $h_g\in H^f_g$ then we set
\begin{equation}
\label{b,,14}
\pi(a)h_g =
\begin{cases}
ah_g, & \ \ g\ge 0,\\
\al_{-g}(a)h_g, &\ \ g\le 0.
\end{cases}
\end{equation}
The scheme of the action of the operator  $\pi(a)$ is presented in
the following  diagram.
\\
\setlength{\unitlength}{0.95mm}
\centerline{
\begin{picture}(80,27)(-37,-10)
\qbezier[300](-80,0)(-15,0)(62,0)\put(65,-1.5){$\bigoplus_{g\in \Gamma} H^f_g $}
\put(70,8){$\pi(a)$}
\put(-50,0){\circle*{1}}\put(-52,-5){$H_{-2x}^f$}
\put(-30,0){\circle*{1}}\put(-32,-5){$H_{-x}^f$}
\put(-10,0){\circle*{1}}\put(-12,-5){$H_0^f$}
\put(10,0){\circle*{1}}\put(8,-5){$H_{x}^f$}
\put(30,0){\circle*{1}}\put(28,-5){$H_{2x}^f$}
\put(47,7){$\dots$}
\put(29,7){$a\cdot$}
\put(9,7){$a\cdot$}
\put(-11,7){$a\cdot$}
\put(-36,7){$\al_{x}(a)\cdot$}
\put(-57.5,7){$\al_{2x}(a)\cdot$}
\put(-77,7){$\dots$}
\end{picture}}
Let us verify  equalities  (\ref{b,,4}) for the representation $\pi$.
Take any  $x\in \Gamma^+$ and  $h_g \in H^f_g$.
\\
If  $x\leq g$ then
\begin{gather*}
\pi(\al_x(a))h_g =\al_{x}(a)h_g,\\[6pt]
U_x\pi(a)U_x^*h_g =\al_x(a L_x(h_g)) =\al_x(a)h_g\al_x(1),
\end{gather*}
and moreover (\ref{b,,8}), inequality $x\leq g$ and the definition
of a transfer operator   imply that  the element
$\al_x(a)h_g\al_x(1)$ coincides with $\al_x(a)h_g$ in the space
$H^f_g$.
\\
For  $0\leq g\leq x$ we have
\begin{gather*}
\pi(\al_x(a))h_g =\al_{x}(a)h_g,\\[6pt]
U_x\pi(a)U_x^*h_g =\al_{x}(1)\al_g(\al_{x-g}(a)L_g(h_g) )=\al_x(a)\al_g(1)h_g\al_g(1)=\al_x(a)h_g\al_g(1),
\end{gather*}
where we used the inequality $\al_x(1)\leq \al_g(1)$. The same argument as  above shows that $\al_x(a)h_g\al_g(1)$
coincides with $\al_x(a)h_g$ in the
space $H^f_g$.
\\
For  $ g\leq 0$ we have
\begin{gather*}
\pi(\al_x(a))h_g =\al_{x-g}(a)h_g,\\[6pt]
U_x\pi(a)U^*_xh_g =\al_{x-g}(1)\al_{x-g}(a)\al_{x-g}(1)h_g=\al_{x-g}(a)h_g.
\end{gather*}
Thus we have proved that $U_x\pi(a)U_x^* =\pi(\al_x(a))$ for any $a\in\cal A$.
It can be shown similarly   that $U_x^*\pi(a)U_x =\pi(L_x(a))$. Indeed, if $0\le g$ one  then  has
\begin{gather*}
\pi(L_x(a))h_g =L_x(a)h_g,\\[6pt]
U^*_x\pi(a)U_xh_g =L_x(a\al_{x}(h_g))=L_x(a)h_g.
\end{gather*}
For  $-x\leq g\leq 0$ one has
\begin{gather*}
\pi(L_x(a))h_g =\al_{-g}(L_x(a))h_g=\al_{-g}(L_{-g}(L_{x+g}(a)))h_g=\al_{-g}(1)L_{x+g}(a)\al_{-g}(1)h_g,\\[6pt]
U^*_x\pi(a)U_xh_g =L_{g+x}(a\al_{x}(1)\al_{g+x}(h_g))=L_{g+x}(a\al_{x+g-g}(1))(h_g)=L_{x+g}(a)\al_{-g}(1)h_g,
\end{gather*}
and moreover (\ref{b,,9}) implies that $\al_{-g}(1)L_{x+g}(a)\al_{-g}(1)h_g$ coincides with $L_{x+g}(a)\al_{-g}(1)h_g$ in the
space $H^f_g$.\\
For  $g\leq -x$ we have
\begin{gather*}
U^*_x\pi(a)U_xh_g =\al_{-g}(1)\al_{-g-x}(a)\al_{-g}(1)h_g,\\[6pt]
\pi(L_x(a))h_g =\al_{-g}(L_x(a))h_g=\al_{-g-x}(\al_{x}(L_{x}(a)))h_g=\al_{-g}(1)\al_{-g-x}(a)\al_{-g}(1)h_g
\end{gather*}
where we used the inequality $0 \leq -g - x$.

To finish the proof it is enough to observe the faithfulness of the
representation $\pi$. But this follows from the definition of the
inner product in (\ref{b,,7}),
the definition of $\pi$ (see the second
line in the diagram) and the standard Gelfand-Naimark faithful
representation of a $C^*$-algebra.
The proof is complete.
\qed
\section{The Banach $^*$-algebra $l_1(\Gamma,\al,\A)$}\label{l1}

This is a starting section for the construction of the principal
object of the article --- the crossed product (that will be
introduced in Section \ref{Cr}). The algebra presented  here
serves as an explicitly described base  of the crossed product.
Hereafter we presume that $(\A,\Gamma^+,\al)$ is a finely
representable $C^*$-dynamical system  and  $L$ is the unique
transfer action (this $L$ does exist by Theorem \ref{b..1} and is
unique by Theorem \ref{complete}).

\begin{I}
\label{l_1 algebra}
Let  $l_1(\Gamma,\al,\A)$ be the set  consisting  of the elements of the form $ a=\{a_g\}_{g\in  \Gamma}$ where
\begin{equation}\label{the choice of a_g}
  a_x\in \A\al_x(1)\quad \textrm{ and  }\quad a_{-x}\in \al_x(1)\A,\qquad \textrm{ for } x\geq 0,
\end{equation}
and such that $\sum_{g\in \Gamma} \|a_g\|< \infty $.
We define the
addition, multiplication by scalar and
involution on $l_1(\Gamma,\al,\A)$ in an obvious manner. Namely, let  $a=\{a_g\}_{g\in  \Gamma}$,
$b=\{b_x\}_{x\in G}\in l_1(\Gamma,\al,\A)$, and let $\lambda\in\mathbb{C}$. We set
\begin{equation}\label{add} (a+b)_g:=a_g+b_g, \end{equation}
\begin{equation}\label{mulscal}
 (\lambda a)_g:=\lambda a_g,
\end{equation}
\begin{equation}\label{invol}
 (a^*)_g:=a^*_{-g}.\end{equation}
Clearly,  these operations are well defined and thus we have equipped $l_1(\Gamma,\al,\A)$ with
the structure of Banach space with isometric involution, the norm taken into account is of course the one given by
$\|a\|=\sum_{g\in \Gamma} \|a_g\|$.
\\
Unfortunately the multiplication  of two elements from
$l_1(\Gamma,\al,\A)$ is not so nice.
It generalizes the convolution multiplication in crossed products by automorphisms and by partial automorphisms.
Its 'strange' structure reflects the arising further 'antisymmetry' between the operators
$U_x$ and $U^*_x$ (see Proposition \ref{universality}) which is mainly due to the fact that in view of (\ref{c11})
one can move $U_x$ only to the right while $U^*_x$ can be moved only to the left.

We put
$$\label{mul}
 (a\cdot b)_g:=\begin{cases}
\displaystyle{\sum_{g=x-y \atop x,y>0}}a_x\al_{x-y}(b_{-y})+
\displaystyle{\sum_{g=y-x\atop x,y>0}}L_x(a_{-x}b_y)+\displaystyle{\sum_{g=x+y \atop x,y \geq 0}}a_x\al_x(b_y),&  {\rm if }\, \, 0 \leq g,\\
\displaystyle{\sum_{g=x-y \atop x,y>0 }}\al_{y-x}(a_x)b_{-y}+
\displaystyle{\sum_{g=y-x \atop x,y>0 }}L_y(a_{-x}b_y)+\displaystyle{\sum_{g=-x-y \atop x,y \geq 0}}\al_y(a_{-x})b_{-y},&   {\rm if }\, \, g < 0.
\end{cases}
$$
where $x,y$ in sums run through $\Gamma^{+}$.
\end{I}
 \begin{Pn}
The above  multiplication  is well defined, and
$$
\|a\cdot b\|\leq \|a\| \cdot \|b\|,\qquad \textrm{for all }\quad a,b\in \A.
$$
\end{Pn}
{\bfseries Proof.}
We need to show that $(a\cdot b)$ satisfies  relations \eqref{the choice of a_g}. To this end take $g\in \Gamma$
and assume that  $0\leq g$ (the case $g<0 $ can be considered in a  similar way). \\
Let $x,y\geq 0$. If $g=x-y$ then, since $\al_g$ is a morphism, we
have
$$a_x\al_{x-y}(b_{-y})=a_x\al_{x-y}(b_{-y})\al_{x-y}(1)\in \A\al_g(1).
$$
If $g=y-x$ then, since $b_y=b_y\al_y(1)$ and $L_x$ is a transfer operator , we have
$$
L_x(a_{-x}b_y)=L_x\big(a_{-x}b_y\al_x(\al_{y-x}(1))\big)=L_x(a_{-x}b_y)\al_{y-x}(1)\in \A\al_g(1).
$$
If $g=x+y$ then, since $b_y=b_y\al_y(1)$ and $\al_x$ is a morphism, we have
$$
a_x\al_{x}(b_{y})=a_x\al_{x}(b_{y}\al_y(1))= a_x\al_{x}(b_y)\al_{x+y}(1)\in \A\al_g(1).
$$
Thus $(a\cdot b)_g \in \A\al_x(1)$.\\
Now, using the fact that $\|\al\|=\|L\|\leq 1$ we have
$$
\|a\cdot b\|=\sum_{g\in \Gamma}\|(a\cdot b)_g\|\leq
$$
$$
\sum_{g\in \Gamma}\Big(\displaystyle{\sum_{g=x-y}}\|a_x\|\|b_{-y}\|+
\displaystyle{\sum_{g=y-x}}\|a_{-x}\|\|b_y\|+\displaystyle{\sum_{g=x+y \atop -g=x+y}}\|a_x\|\|b_y\|\Big)=
$$
$$
\Big(\sum_{g\in \Gamma}\|a_g\|\Big) \Big(\sum_{g\in \Gamma}\|b_g\|\Big)=\|a\|\|b\|.
$$
Hence the multiplication is well defined and the proof is complete.
\qed

Let us provide the following notational convention: \\
Let $g_0$ be in $\Gamma$ and let  $a$ be in $\A\al_{g_0}(1)$, if $g_0\geq 0$, or
in $\al_{-g_0}(1)\A $, if $g_0<0 $. We shall denote by $a\delta_{g_0}$ the  element  given by
$$
(a\delta_{g_0})_g=a\delta_{(g_0,g)}
$$
where $\delta_{(g_0, g)}$ is the Kronecker symbol.
\\
Then the elements $\al_x(1)\delta_x$,  $x\in \Gamma^{+}$, and the algebra $\A\delta_0$ generate
(with the help of the above defined operations) a dense subspace (in fact a $^*$-subalgebra) of $l_1(\Gamma,\al,\A)$.
 Furthermore, we have the natural embedding of $\A$ into $l_1(\Gamma,\al,\A)$ given by
$$
\A \ni a \longmapsto a\delta_0.
$$
It is easy to check that  under this embedding  the unit $1\in \A$ coincides with  the element $1\delta_0$ which is
   neutral  with respect to the  multiplication we have defined.
\begin{Tm}\label{algebrawithinvolution}
The set $l_1(\Gamma,\al,\A)$ with the above defined operations  becomes a unital $^*$-Banach  algebra.
\end{Tm}
{\bfseries Proof.} It is enough to  verify   the equality  $(a\cdot b)^*=b^*\cdot a^*$ and the associativity
 of multiplication (the distribution laws  are readily checked because $\al$ and $L$ are linear).
\\
Let us prove the first property. Let $g$ be in $\Gamma$. If $   g\geq 0$, then using the  positivity of $L$ we have
$$
 (a\cdot b)_{-g}^*=\displaystyle{\sum_{-g=x-y}}b_{-y}^*\al_{y-x}(a_x^*)+
 \displaystyle{\sum_{-g=y-x}}L_y(b_y^*a_{-x}^*)+\displaystyle{\sum_{-g=-x-y}}b_{-y}^*\al_y(a_{-x}^*),
$$
and simply by definition
$$
 (b^*\cdot a^*)_{g}=\displaystyle{\sum_{g=x-y}}b_{-x}^*\al_{x-y}(a_{-y}^*)+
 \displaystyle{\sum_{g=y-x}}L_y(b_x^*a_{-y}^*)+\displaystyle{\sum_{g=x+y}}b_{-x}^*\al_y(a_{-y}^*).
$$
Replacing $x$ by $y$, one sees that $(a\cdot b)_{-g}^*=(b^*\cdot a^*)_{g}$.
Using the same argument for $g \leq 0$ one obtains  the desired equality:  $(a\cdot b)^*=b^*\cdot a^*$.
\\
Clearly, to show the associativity it  suffices to consider the elements of the form
$a\delta_{g_1}$, $b\delta_{g_2}$, $c\delta_{g_3} \in l_1(\Gamma,\al,\A)$ where $g_1,g_2,g_3\in \Gamma$ are fixed.
 However, anyway we are   faced against a number of possibilities which have to be checked.
 This may be a source of a pleasure as well as a cause of a headache
  hence we confine  ourselves to the case when $g_1+g_2+g_3\geq 0$ and leave the opposite case to the enthusiasts.
\\
Suppose that $g_1+g_2+g_3\geq 0$. The routine computation shows that:
\\
1) If $g_1,g_2,g_3\geq 0$,   then
$$
(a\delta_{g_1}\cdot b\delta_{g_2})\cdot c\delta_{g_3}= (a\al_{g_1}(b)\delta_{g_1+g_2})\cdot c\delta_{g_3}
= a\al_{g_1}(b)\al_{g_1+g_2}(c)\delta_{g_1+g_2+g_3}
$$
$$
= a\al_{g_1}(b\al_{g_2}(c)) \delta_{ g_1+g_2+g_3}=a\delta_{g_1}\cdot (b\al_{g_2}(c)
\delta_{g_2+g_3})=a\delta_{g_1}\cdot (b\delta_{g_2}\cdot c\delta_{g_3}).
$$
2) If $g_1<0$ and $g_2,g_3\geq 0$,  then
$$
(a\delta_{g_1}\cdot b\delta_{g_2})\cdot c\delta_{g_3}=
\begin{cases}L_{-g_1}(ab)\al_{g_1+g_2}(c)\delta_{g_1+g_2+g_3}, & \textrm{when}\,\, \, g_1+g_2 \geq 0
\\
L_{-g_1-g_2}(L_{g_2}(ab)c)\delta_{g_1+g_2+g_3}, & \textrm{when}\,\, \, g_1+g_2 <0
\end{cases}
$$
$$
=\begin{cases}
L_{-g_1}\big(ab\al_{-g_1}(\al_{g_1+g_2}(c))\big)\delta_{g_1+g_2+g_3}, & \textrm{when}\,\, \, g_1+g_2 \geq 0
\\
L_{-g_1-g_2}\big(L_{g_2}(ab\al_{g_2}(c))\big)\delta_{g_1+g_2+g_3}, & \textrm{when}\,\, \, g_1+g_2 <0
\end{cases}
$$
$$
= L_{-g_1}(ab\al_{g_2}(c))\delta_{g_1+g_2+g_3}=a\delta_{g_1}\cdot (b\delta_{g_2}\cdot c\delta_{g_3})
$$
\\
3) If $g_2<0$ and $g_1,g_3\geq 0$,  then
$$
(a\delta_{g_1}\cdot b\delta_{g_2})\cdot c\delta_{g_3}=
\begin{cases}a\al_{g_1+g_2}(b)\al_{g_1+g_2}(c)\delta_{g_1+g_2+g_3}, & \textrm{when}\,\, \, g_1+g_2 \geq 0
\\
L_{-g_1-g_2}(\al_{-g_1-g_2}(a)bc)\delta_{g_1+g_2+g_3}, & \textrm{when}\,\, \, g_1+g_2 <0
\end{cases}
$$
$$
=
\begin{cases}a\al_{g_1+g_2}(bc)\delta_{g_1+g_2+g_3}, & \textrm{when}\,\, \, g_1+g_2 \geq 0
\\
aL_{-g_1-g_2}(bc)\delta_{g_1+g_2+g_3}, & \textrm{when}\,\, \, g_1+g_2 <0
\end{cases}
$$
and
$$
a\delta_{g_1}\cdot (b\delta_{g_2}\cdot c\delta_{g_3})=
\begin{cases} a\al_{g_1}(L_{-g_2}(bc))\delta_{g_1+g_2+g_3}, & \textrm{when}\,\, \, g_2+g_3 \geq0
\\
a\al_{g_1+g_2+g_3}(L_{g_3}(bc))\delta_{g_1+g_2+g_3}, & \textrm{when}\,\, \, g_2+g_3 < 0
\end{cases}.
$$
We have  the four (sub)possibilities:
\par
3a)   For $g_1+g_2 \geq 0 $ and $g_2+g_3 < 0$ we have
$$
\big(a\delta_{g_1}\cdot (b\delta_{g_2}\cdot c\delta_{g_3})\big)_{g_1+g_2+g_3}=
a\al_{g_1+g_2}\big(\al_{g_3}(L_{g_3}(bc)\big)=a\al_{g_1+g_2}\big(\al_{g_3}(1)bc\al_{g_3}(1)\big)
$$
$$
=a\al_{g_1+g_2+g_3}(1) \al_{g_1+g_2}(bc)=a\al_{g_1+g_2}(bc)=\big((a\delta_{g_1}\cdot b\delta_{g_2})\cdot c\delta_{g_3}\big)_{g_1+g_2+g_3}
$$
\par
3b)   For $g_1+g_2 \geq 0 $ and $g_2+g_3 \geq 0$ we have
$$
\big(a\delta_{g_1}\cdot (b\delta_{g_2}\cdot c\delta_{g_3})\big)_{g_1+g_2+g_3}=a\al_{g_1+g_2}
\big(\al_{-g_2}(L_{-g_2}(bc)\big)=a\al_{g_1+g_2}\big(\al_{-g_2}(1)bc\al_{-g_2}(1)\big)
$$
$$
=a\al_{g_1+g_2}(bc\big)=\big((a\delta_{g_1}\cdot b\delta_{g_2})\cdot c\delta_{g_3}\big)_{g_1+g_2+g_3}
$$
\par
3c)   For $g_1+g_2 < 0 $ and $g_2+g_3 \geq 0$ we have
$$
\big(a\delta_{g_1}\cdot (b\delta_{g_2}\cdot c\delta_{g_3})\big)_{g_1+g_2+g_3}=a\al_{g_1}\big(L_{g_1}(L_{-g_1-g_2}(bc)\big)=
a\al_{g_1}(1)L_{-g_1-g_2}(bc)\al_{g_1}(1)
$$
$$
=aL_{-g_1-g_2}(bc\al_{-g_2}(1))=aL_{-g_1-g_2}(bc)=\big((a\delta_{g_1}\cdot b\delta_{g_2})\cdot c\delta_{g_3}\big)_{g_1+g_2+g_3}
$$
\par
3d)   For $g_1+g_2 < 0 $ and $g_2+g_3 < 0$ we have
$$
\big(a\delta_{g_1}\cdot (b\delta_{g_2}\cdot c\delta_{g_3})\big)_{g_1+g_2+g_3}=a\al_{g_1+g_2+g_3}\big(L_{g_1+g_2+g_3}(L_{-g_1-g_2}(bc))\big)
$$
$$
=a\al_{g_1+g_2+g_3}(1)L_{-g_1-g_2}(bc)\al_{g_1+g_2+g_3}(1)=
aL_{-g_1-g_2}(bc\al_{g_3}(1))=\big((a\delta_{g_1}\cdot b\delta_{g_2})\cdot c\delta_{g_3}\big)_{g_1+g_2+g_3}
$$
\\
4) If $g_3<0$ and $g_1,g_2\geq 0$,  then
$$
(a\delta_{g_1}\cdot b\delta_{g_2})\cdot c\delta_{g_3}=
\begin{cases}a\al_{g_1}(b\al_{g_2+g_3}(c))\delta_{ g_1+g_2+g_3}, & \textrm{when}\,\, \, g_2+g_3 \geq 0
\\
a\al_{g_1+g_2+g_3}(\al_{-g_2-g_3}(b)c)\delta_{g_1+g_2+g_3}, & \textrm{when}\,\, \, g_2+g_3 <0
\end{cases}
$$
$$
= a\al_{g_1}(b)\al_{g_1+g_2+g_3}(c))\delta_{g_1+g_2+g_3}=a\delta_{g_1}\cdot (b\delta_{g_2}\cdot c\delta_{g_3})
$$
\\
5) If $g_1,g_2 <0$ and $g_3\geq 0$,  then as $g_1+g_2+g_3\geq 0$  we have $g_2+g_3>0$, and thus
$$
a\delta_{g_1}\cdot (b\delta_{g_2}\cdot c\delta_{g_3})=
L_{-g_1}(aL_{-g_2}(bc))\delta_{g_1+g_2+g_3}=
L_{-g_1}(L_{-g_2}(\al_{-g_2}(a)bc))\delta_{g_1+g_2+g_3}
$$
$$
=L_{-g_1-g_2}(\al_{-g_2}(a)bc)\delta_{g_1+g_2+g_3}=(a\delta_{g_1}\cdot b\delta_{g_2})\cdot c\delta_{g_3}
$$
\\
6) If $g_1,g_3 < 0$ and $g_2\geq 0$,  then as $g_1+g_2+g_3\geq 0$  we have $g_2+g_3>0$ and $g_1+g_2>0$. Thus
$$
a\delta_{g_1}\cdot (b\delta_{g_2}\cdot c\delta_{g_3})=L_{-g_1}(ab \al_{g_2+g_3}(c))\delta_{g_1+g_2+g_3}
$$
$$
=L_{-g_1}\big(ab \al_{-g_1}(\al_{g_1+g_2+g_3}(c))\big)\delta_{g_1+g_2+g_3}=L_{-g_1}(ab)\al_{g_1+g_2+g_3}(c)\delta_{g_1+g_2+g_3}
$$
$$
=(a\delta_{g_1}\cdot b\delta_{g_2})\cdot c\delta_{g_3}.
$$
\\
7) If $g_2,g_3 < 0$ and $g_1\geq 0$,  then as $g_1+g_2+g_3\geq 0$  we have $g_1+g_2>0$, and thus
$$
a\delta_{g_1}\cdot (b\delta_{g_2}\cdot c\delta_{g_3})=
a\al_{g_1+g_2+g_3}(\al_{-g_3}(b)c)\delta_{g_1+g_2+g_3}
$$
$$
=a\al_{g_1+g_2}(b)\al_{g_1+g_2+g_3}(c)\delta_{ g_1+g_2+g_3}=(a\delta_{g_1}\cdot b\delta_{g_2})\cdot c\delta_{g_3}.
$$
Thus the equality $a\delta_{g_1}\cdot (b\delta_{g_2}\cdot c\delta_{g_3})=
(a\delta_{g_1}\cdot b\delta_{g_2})\cdot c\delta_{g_3}$, in the case $g_1+g_2+g_3\geq 0$, is proved.
\qed
\section{The crossed product $\A \rtimes_\al \Gamma$}\label{Cr}
In this section we discuss the main object of the paper the crossed product of a finely representable $C^*$-dynamical system.

\begin{Dn}
The crossed product of a finely representable  $C^*$-dynamical system $(\A,\Gamma^{+},\al)$ (see \ref{fine})
is the $C^*$-algebra $\A \rtimes_\al \Gamma$ obtained by taking the enveloping $C^*$-algebra  of $l_1(\Gamma,\al,\A)$ (see \ref{l_1 algebra}).
\end{Dn}%\section{Representations of the crossed-product}
We  aim at the investigation  of the structure of  $\A \rtimes_\al \Gamma$, but before that let us  justify the
above definition and show  that $\A \rtimes_\al \Gamma$ is universal with respect to covariant representations.
\begin{Dn}\label{covrep}
Let $(\A,\Gamma^+, \al)$ be finely representable. A triple $(\pi,H,U)$ consisting of a Hilbert space~$H$, a
non-degenerate representation $\pi\!:{\cal A}\to L(H)$ and a semigroup homomorphism $U\!:\Gamma^+ \to L(H)$,
is a \emph{covariant representation} of $(\A,\Gamma^+, \al)$, if
 for every $a\in \A$ and $x\in\Gamma^+$ we have
$$
U_x\pi(a)U^*_x=\pi(\alpha_x(a)),
\qquad  U^*_x\pi(a)U_x=\pi(L_x(a)).
$$
\end{Dn}
\begin{Pn}\label{universality}
Let $(\pi,U,H)$ be a covariant representation of $(\A,\Gamma^{+},\al)$. Then the formula
$$
(\pi\times U)\Big(\sum_{x>0} a_{-x}\delta_{-x} +a_0\delta_{0}+  \sum_{x>0} a_x\delta_{x}\Big):=\sum_{x>0}U_x^* a_{-x} +a_0+  \sum_{x>0} a_xU_{x},
$$
 defines a representation of $l_1(\Gamma,\al,\A)$ and hence establishes a representation of $\A \rtimes_\al \Gamma$.
\end{Pn}
{\bfseries Proof.} Clearly, $(\pi\times U)$ is   linear and
preserves  the involution. In order to show that $(\pi\times U)$
is  multiplicative  let us consider the elements $a\delta_{g_1}$,
$b\delta_{g_2}\in l_1(\Gamma,\al,\A)$ (we recall that if
$a\delta_g\in l_1(\Gamma,\al,\A)$ then  $a\in  \A\al_g(1)$ when
$g\leq 0$ and $a\in\al_{-g}(1)\A$ otherwise). We have the
following possibilities:
\\
\ \ I) $g_1+g_2\geq 0$, in other words  $g_1\geq-g_2$,\ \   $g_2\geq -g_1$.
\\

1)
If $g_1> 0$,  $g_2< 0$ then
$$
(\pi\times U)(a\delta_{g_1})(\pi\times U)(b\delta_{g_2})=\pi(a)U_{g_1}U_{-g_2}^*\pi(b)=\pi(a)U_{g_1+g_2}U_{-g_2}U_{-g_2}^*\pi(b)
$$
$$
=\pi(a)U_{g_1+g_2}\pi(\al_{-g_2}(1)b)=\pi(a)U_{g_1+g_2}\pi(b)U_{g_1+g_2}^*U_{g_1+g_2}
$$
$$
=(\pi\times U)(a\al_{g_1+g_2}(b)\delta_{g_1+g_2})=(\pi\times U)\big((a\delta_{g_1})(b\delta_{g_2})\big).
$$
\par
2)
If $g_1<0$, $g_2> 0$ then
$$
(\pi\times U)(a\delta_{g_1})(\pi\times U)(b\delta_{g_2})=U_{-g_1}^*\pi(a)\pi(b)U_{g_2}=U_{-g_1}^*\pi(ab)U_{-g_1}U_{g_1+ g_2}
$$
$$
=\pi(L_{-g_1}(ab))U_{g_1+g_2}=(\pi\times U)(L_{-g_1}(ab)\delta_{g_1+g_2})=(\pi\times U)\big((a\delta_{g_1})(b\delta_{g_2})\big).
$$
\par
3) If $g_1,g_2\geq 0$ then
$$
(\pi\times U)(a\delta_{g_1})(\pi\times U)(b\delta_{g_2})=\pi(a)U_{g_1}\pi(b)U_{g_2}=\pi(a)U_{g_1}\pi(b)U^{*}_{g_1}U_{g_1}U_{g_2}
$$
$$
=\pi(a)\pi(\al_{g_1}(b))U_{g_1+g_2}=(\pi\times U)(a\al_{g_1}(b)\delta_{g_1+g_2})=(\pi\times U)\big((a\delta_{g_1})(b\delta_{g_2})\big).
$$
\\
\ \ II) $g_1+g_2< 0$, in other words  $g_1<-g_2$,\ \    $g_2< -g_1$.
\\

1)
If $g_1> 0$,  $g_2< 0$ then
$$
(\pi\times U)(a\delta_{g_1})(\pi\times U)(b\delta_{g_2})=\pi(a)U_{g_1}U_{-g_2}^*\pi(b)=\pi(a)U_{g_1}U_{g_1}^*U_{-g_2-g_1}^*\pi(b)
$$
$$
=\pi(a\al_{g_1}(1))U_{-g_1-g_2}^*\pi(b)=U_{-g_1-g_2}^*U_{-g_1-g_2}\pi(a)U_{-g_1-g_2}^*\pi(b)
$$
$$
=(\pi\times U)(\al_{-g_1-g_2}(a)b\delta_{g_1+g_2})=(\pi\times U)\big((a\delta_{g_1})(b\delta_{g_2})\big).
$$
\par
2)
If $g_1<0$, $g_2> 0$ then
$$
(\pi\times U)(a\delta_{g_1})(\pi\times U)(b\delta_{g_2})=U_{-g_1}^*\pi(a)\pi(b)U_{g_2}=U_{-g_1-g_2}^*U_{g_2}^*\pi(ab)U_{g_2}
$$
$$
=U_{-g_1-g_2}\pi(L_{g_2}(ab))=(\pi\times U)(L_{g_2}(ab)\delta_{g_1+g_2})=(\pi\times U)\big((a\delta_{g_1})(b\delta_{g_2})\big).
$$
\par
3) If $g_1,g_2< 0$ then
$$
(\pi\times U)(a\delta_{g_1})(\pi\times U)(b\delta_{g_2})=U_{-g_1}^*\pi(a)U_{-g_2}^*\pi(b)=U_{-g_1}^*U_{-g_2}^*U_{-g_2}\pi(a)U_{-g_2}^*\pi(b)
$$
$$
=U_{-g_1-g_2}^*\pi(\al_{-g_2}(a))\pi(b)=(\pi\times U)(a\al_{g_1}(b)\delta_{g_1+g_2})=(\pi\times U)\big((a\delta_{g_1})(b\delta_{g_2})\big).
$$
\qed

As we know that there exists a  covariant representation $(\pi,U,H)$ such that $\pi$ is
 faithful (see Theorem \ref{b..1}) this implies the existence of  a representation of $l_1(\Gamma,\al,\A)$ which is faithful on
 $\A\delta_0$ and hence the algebra $\A$ is naturally embedded into the crossed product  $\A \rtimes_\al \Gamma$.
 The argument we have just used has in fact much stronger consequences - see  item $(iv)$ of the following
\begin{Tm}\label{crossedstructure}
Let  $u_x$ (resp. $u_{-x}$) be  the element of  $\A \rtimes_\al \Gamma$ corresponding
 to the element $\al_x(1)\delta_x $ (resp. $\al_x(1)\delta_{-x} $)  of $l_1(\Gamma,\al,\A)$, $x\in \Gamma^{+}$. Then
\begin{itemize}
\item[(i)] the family $\{u_x\}_{x\in\Gamma^{+}}$ forms a semigroup of partial isometries,
\item[(ii)] for each $x\in \Gamma^{+}$ and $a\in \A$ we have
$$
u_x^*=u_{-x},\qquad u_xa u_x^*=\al_x(a),\quad u_x^*a u_x=L_x(a).
$$
\item[(iii)] the elements of the form
\begin{equation}\label{suma}
a=\sum_{x\in F} u_x^*a_{-x} +a_0+  \sum_{x\in F} a_xu_x , \quad \ F\subset \Gamma^+\setminus\{0\}, \ \ \vert F\vert < \infty
\end{equation}
where $a_{x}\in \A u_xu^*_x$, $a_{-x}\in u_xu^*_x\A$, form a dense $^*$-subalgebra $ C_0$ of $\A \rtimes_\al \Gamma$.
\item[(iv)] if $a$ is of the form \eqref{suma}, then  we have the inequalities
$$
\| a_g \|\leq \|a\|, \qquad g\in \Gamma.
$$
In particular the  'coefficients' $a_g$, $g\in \Gamma$,  of $a$ are determined in a unique way.
\end{itemize}
\end{Tm}
{\bfseries Proof.} For $x,y\in \Gamma^{+}$ we have $\al_x(1)\delta_x \al_y(1)\delta_y=\al_x(1)\al_x(\al_y(1))\delta_{x+y}=\al_{x+y}(1)\delta_{x+y}$,
whence the family $\{u_x\}_{x\in\Gamma^{+}}$ forms a semigroup, that is $(i)$ holds. \\
As $(\al_x(1)\delta_x)^*=\al_x(1)\delta_{-x}$ we have $
u_x^*=u_{-x}$. For $a\in \A$ we have
$$
(\al_x(1)\delta_x)( a\delta_0) (\al_x(1)\delta_{-x})=(\al_x(a)\delta_x ) (\al_x(1)\delta_{-x})=\al_x(a)\delta_0,
$$
$$
(\al_x(1)\delta_{-x}) (a\delta_0) (\al_x(1)\delta_{x})=(\al_x(1) a\delta_{-x}) (\al_x(1)\delta_{x})
=L_x(\al_x(1) a \al_x(1))\delta_{0}= L_x(a)\delta_{0},
$$
which implies  $(ii)$, and in particular it follows that $u_x$  is a partial isometry.
\\
As element $\sum_{x\in F} u_x^*a_{-x} +a_0+  \sum_{x\in F} a_xu_x$ corresponds to the element
$\sum_{x\in F} a_{-x}\delta_{-x} +a_0\delta_{0}+  \sum_{x\in F} a_x\delta_{x}$ in $l_1(\Gamma,\al,\A)$
it follows that  $(iii)$ is true.
\\
 Now let us verify  $(iv)$ for $g=0$. To this end, take any  covariant representation $(\pi,U,H)$ such that $\pi$ is  faithful.
Consider the space ${\cal H} =l_2 (\Gamma, H)$ and the
representation  $\nu \! : \A \rtimes_\al \Gamma\to L({\cal
H})$  given by the formulae
\begin{gather*}
(\nu (a)\xi )_g = \pi (a) (\xi_g), \qquad\text{where}\quad
a\in {\cal A}, \quad l_2 (\Gamma, H) \ni \xi = \{ \xi_g  \}_{g\in
\Gamma}\,;\\[6pt]
(\nu (u_x)\xi )_g = U_x (\xi_{g-x}),\qquad
(\nu (u^*_x)\xi )_g= U^*_x (\xi_{g+x}).
\end{gather*}
Routine verification shows that $\nu (\A)$ and $\nu
(u_x)$ satisfy all the conditions of a covariant representation and thus by
  Proposition \ref{universality} $\nu$ is well defined.
\\
Now take any $a\in  \A \rtimes_\al \Gamma$ given by
\eqref{suma} and for a given $\varepsilon > 0$ chose  a vector
$\eta \in H$ such that
\begin{equation}\label{e}
\Vert \eta \Vert =1 \quad \text{and}\quad \Vert \pi (a_0)
\eta \Vert > \Vert \pi (a_0)  \Vert - \varepsilon .
\end{equation}
Set $\xi \in l_2 (\Gamma, H)$ by $\xi_g = \delta_{(0,g)}\eta $,
where $\delta_{(i,j)}$ is the Kronecker symbol. We have that $\Vert
\xi \Vert = 1$ and the explicit form of $\nu (a) \xi$  and
\eqref{e} imply
\[
\Vert \nu (a) \xi \Vert \ge \Vert \pi (a_0) \eta \Vert >
\Vert \pi (a_0)  \Vert - \varepsilon
\]
which by the arbitrariness of $\varepsilon$ proves the desired
inequality: $
\Vert a  \Vert \ge \Vert \nu(a)  \Vert\ge \Vert \pi (a_0)  \Vert = \Vert  a_0 \Vert.
$
\\
In order to verify the corresponding inequality for an arbitrary $g\in \Gamma$ take any $x\in\Gamma^+$ and observe that
$$
    (u_xa)_0= u_xu^*_xa_{- x}=a_{-x}, \qquad (a u_x^*)_0=  a_{ x}u_xu^*_x=a_x.
 $$
 Hence, we have
 $$\|a_{-x}\|\leq \|u_xa\|\leq \|a\|,\qquad \|a_{x}\|\leq \|a u_x^*\|\leq \|a\|.$$
  \qed
\begin{Cy}
We have the one-to-one correspondence:
$$
(\pi,U,H)\longleftrightarrow (\pi\times U, H)
$$
between covariant representations of $(\A,\Gamma^+,\al)$ and non-degenerate representations of $\A \rtimes_\al \Gamma$.
\end{Cy}
{\bfseries Proof.} It follows from Proposition \ref{universality} and items $(i)$, $(ii)$ of the above theorem.
\qed
\begin{I}\label{conditionexpect}
Items $(iii)$ and $(iv)$ of Theorem \ref{crossedstructure}  mean that  one
can define the linear and continuous maps $E_x \!: C_0 \to
\A u_x u^*_x$ and $E_{-x} \!: C_0 \to u_x u_x^*\A$, \,
$x\in\Gamma^+$ given by
$$
E_g(a)= a_g, \qquad \quad a\in C_0,\,\, g\in \Gamma.
$$
By continuity these mappings can be expanded onto the whole of
$\A \rtimes_\al \Gamma$ thus defining the
{\bfseries 'coefficients'} of an arbitrary element $a\in \A \rtimes_\al \Gamma$. We shall show further in Theorem \ref{uniqueNk}
  that these
coefficients determine $a$ in a unique way.
\end{I}
 The next theorem  shows that  the norm of an
element $a\in C_0$ can be calculated only in terms of the elements
of $\A$ (0-degree coefficients of the powers of $aa^*$).
 \begin{Tm}
\label{3a.N} Let $a\in C_0\subset \A \rtimes_\al \Gamma$ be of the form \eqref{suma}. Then we have
\begin{equation}
\label{be3.131} \Vert a \Vert = \lim_{k\to\infty} \sqrt[4k]{ \left\Vert E_0 \left[
(a\cdot a^*)^{2k}\right]\right\Vert }.
\end{equation}
\end{Tm}
{\bfseries Proof.}   Applying the known equality $\left\Vert  \sum_{i=1}^m d_i  \right\Vert^2 \le m \left\Vert  \sum_{i=1}^m d_i
d_i^* \right\Vert$ which holds for any elements $d_1, ..., d_m$ in an arbitrary $C^*$-algebra,
 we obtain
 $$
 \Vert a \Vert^2 \leq
(2|F|+1) \left\Vert \sum_{x\in F} a_{x}u_{x}u^{*}_{x}a_{x}^* + a_0a_0^*+ \sum_{x\in F}u^{*}_{x}a_{-x}a_{-x}^* u_{x} \right\Vert =
(2|F|+1) \Vert E_0 (aa^*)\Vert. $$
  On the other hand as
$E_0$ is contractive we have
$$
 \Vert a \Vert^2 = \Vert aa^* \Vert \ge \Vert E_0
(aa^*)\Vert
$$
 thus
\begin{equation}
\label{be3.101} \Vert E_0 (aa^*)\Vert \le  \Vert aa^* \Vert = \Vert a \Vert^2  \le (2|F|+1) \Vert
E_0 (aa^*)\Vert
\end{equation}
Applying (\ref{be3.101}) to $(aa^* )^k$ and having in mind that $(aa^* )^k = (aa^* )^{k*} $ and
$\Vert (aa^* )^{2k}\Vert = \Vert a \Vert^{4k} $ one has
$$
\Vert E_0 \left[
(aa^*)^{2k}\right]\Vert \le
 \Vert (aa^*)^k \cdot  (aa^*)^{k*}  \Vert =
\Vert a \Vert^{4k}  \le (2|F^k|+1) \Vert E_0 \left[ (aa^*)^{2k}\right]\Vert
$$
where $F^k$ is set of all elements of $\Gamma$ representable as a product of $k$ elements from $F$. We recall that the so called
subexponential groups include the commutative groups and thus $\lim_{k\to\infty} |F^k|^{\frac{1}{k}}=1$.\\ So $$
 \sqrt[4k]{ \left\Vert E_0 \left[ (aa^*)^{2k}\right]\right\Vert }\le \Vert a \Vert \le
\sqrt[4k]{ 2|F^k|+1 }\cdot \sqrt[4k]{  \left\Vert E_0 \left[ (aa^*)^{2k}\right]\right\Vert }
 $$
Observing the equality $$
 \lim_{k\to\infty}\sqrt[4k]{ 2|F^k|+1 } =1
 $$
  we conclude that
  $$
   \Vert a
\Vert = \lim_{k\to\infty} \sqrt[4k]{ \left\Vert E_0 \left[ (aa^*)^{2k}\right]\right\Vert }.
 $$
  The
proof is complete.
\qed

\section{Faithful and regular representations of \\the crossed product}\label{Faith}

Here we give a criterion for a representation of  $\A \rtimes_\al \Gamma$ to be faithful and present the regular representation of
$\A \rtimes_\al \Gamma$. In order to study faithful representations of $\A \rtimes_\al \Gamma$ we introduce the following
\begin{Dn} Let $(\pi ,U, H)$ be a covariant
representation of $(\A,\Gamma^*,\al)$. We shall say that $(\pi ,U, H)$ possesses
{\em property } $(^*)$ if  for any element $a\in C_0$ of the form  (\ref{suma})  we have
$$
\|E_0(a)\| \leq \|(\pi\times U)(a)\|\qquad (^*)
$$
in other words
$
\Vert a _0\Vert \leq \Vert \sum_{x\in F} U_x^*\pi (a_{-x})+\pi(a_0)+ \sum_{x\in F} \pi (a_x)U_x \Vert  .
$
\end{Dn}

Let us observe that if $(\pi ,U, H)$ possess
property  $(^*)$ then the representation $\pi$ is faithful, and the  mapping $
\Nc((\pi\times U)(a)):=E_0(a),$  $a\in C_0$,
 extends uniquely  up to the positive, contractive, conditional expectation from $(\pi\times U)(\A \rtimes_\al \Gamma)$ onto $\A$.
\begin{Tm}\label{isomorphism theorem}
Let $(\pi ,U, H)$ be a covariant
representation of $(\A,\Gamma^+,\al)$. The representation $(\pi\times U)$ of $\A \rtimes_\al \Gamma$ is faithful
iff $(\pi ,U, H)$ possesses property $(^*)$.
\end{Tm}
{\bfseries Proof.} Necessity follows from item $(iv)$ of Theorem \ref{crossedstructure}. Let us show the sufficiency.\\
Take any $a\in C_0$. By Theorem \ref{3a.N} and the definition of property $(^*)$ we have
$$
\label{be3.013} \Vert a \Vert = \lim_{k\to\infty} \sqrt[4k]{ \left\Vert E_0 \left[
(aa^*)^{2k}\right]\right\Vert } = \lim_{k\to\infty} \sqrt[4k]{ \left\Vert \Nc\left[(\pi\times U)
(aa^*)^{2k}\right]\right\Vert }
$$
$$
\leq \lim_{k\to\infty} \sqrt[4k]{ \left\Vert (\pi\times U)
(aa^*)^{2k}\right\Vert }=\lim_{k\to\infty} \sqrt[4k]{ \left\Vert (\pi\times U)
(aa^*)^{k}(\pi\times U)
(aa^*)^{k}\right\Vert }
$$
$$
=\lim_{k\to\infty} \sqrt[4k]{ \left\Vert (\pi\times U)
(a)\right\Vert^{4k} }=\left\Vert (\pi\times U)
(a)\right\Vert.
$$
Hence $
\|a \| =\|(\pi\times U) (a) \| $ on a dense subset of  $\A \rtimes_\al \Gamma$.
\qed
\begin{Cy}\label{dual action}
We have the action of the dual group $\hat{\Gamma}$ on $\A \rtimes_\al \Gamma$ by the automorphisms  given by
$$
\lambda a:=a, \quad a\in \A, \qquad  \lambda u_g :=\lambda_g u_g,\qquad g\in \Gamma,\, \lambda\in \hat{\Gamma},\, \lambda_g=\lambda(g)
$$
(here we consider $\Gamma$ as a discrete group).
\end{Cy}
{\bfseries Proof.}
Suppose that  $\A \rtimes_\al \Gamma$ is faithfully and nondegenerately represented on a Hilbert space $H$.
Then for each $\lambda \in\hat{ \Gamma}$ the triple $(id, \lambda u , H)$ is a covariant representation possessing property $(^*)$,
whence $(id\times \lambda u)$ is an automorphism of $\A \rtimes_\al \Gamma$.

\qed

The next result shows that
 any element $a\in \A \rtimes_\al \Gamma$ can be 'restored' by its
 coefficients $E_g (a)$, $g\in \Gamma$, see  \ref{conditionexpect}.

\begin{Tm}
\label{uniqueNk} Let\/  $a\in \A\times_{\delta}\Gamma$.
Then the following conditions are equivalent:

\smallskip
\quad\llap{$(i)$}\ \ $a=0;$

\smallskip
\quad\llap{$(ii)$}\ \ $E_g (a)=0$, \,$g\in\Gamma;$

\smallskip
\quad\llap{$(iii)$}\ \ $E_0 (a^*a)=0$.

\end{Tm}
 {\bfseries Proof.} Clearly, $(i)$ implies $(ii)$ and  $(iii)$.
 \\
 We now prove $(ii)\Rightarrow (i)$.
Let us suppose that  $\A \rtimes_\al \Gamma$ is faithfully and
nondegenerately represented on a Hilbert space $H$. Thus to prove
that $a=0$ it is enough to show that for any fixed $\xi,\eta \in
H$ with $\Vert \xi \Vert = \Vert \eta \Vert = 1$ we have
\begin{equation}
\label{zero}
\langle a\xi , \eta \rangle = 0
\end{equation}
where $\langle \, ,\, \rangle $ is the inner product in $H$.
\\
Recall that $C_0$ is dense in $\A \rtimes_\al \Gamma$ and hence we can choose a sequence  $a_n $,
$n=1,2, \ldots $  of elements of $C_0$ tending to $a$:
$$
a_n=\sum_{x\in F(n)} u_x^*a_{-x}^{(n)} +a_0^{(n)}+  \sum_{x\in F(n)} a_x^{(n)}u_x ,
$$
where $F(n)$ is finite subset of $\Gamma^+\setminus\{0\}$ and $a_{-x}^{(n)}\in \al_x(1)\A$, $a_{x}^{(n)}\in \A\al_x(1)$.
\\
Let $\lambda \in \hat{\Gamma}$ and consider the elements $\lambda a_n $ (see Corollary \ref{dual action}).
 We  define the sequence
$f_n $, $n=1,2, \ldots,$
of  functions on $\hat{\Gamma}$  by
\begin{equation}
\label{flambda}
f_n (\lambda ) = \langle  \lambda a_n \xi, \eta \rangle  =\sum_{x\in F(n)} \gamma_{-x}^{(n)}\lambda_{-x}
 +\sum_{x\in F(n)\cup\{0\}} \gamma_{x}^{(n)}\lambda_{x}
\end{equation}
where $
\gamma_{-x}^{(n)} = \langle U^{*}_{x}a^{(n)}_{-x}\xi ,\eta \rangle,$ $\gamma_x^{(n)} = \langle a_{x}^{(n)}U_{x} \xi, \eta\rangle $,
 $x\in F(n)\cup\{0\}$. It follows that $f_n$, $n=1,2, \ldots,$  are continuous on $\hat{\Gamma}$
(since $\Gamma$ is discrete $\hat{\Gamma}$ is compact).
\\
%Since $a_n \stackrel{\longrightarrow}{_{_{n\rightarrow \infty}}} a$ it follows (in view of Corollary \ref{dualaction}) that $\lambda a_n \stackrel{\longrightarrow}{_{_{n\rightarrow \infty}}} \lambda a$
%\begin{equation}
%\label{xn>}
%\Vert a_{n_1} (\lambda )- a_{n_2}(\lambda )\Vert =
%\Vert a_{n_1} - a_{n_2}\Vert \stackrel{\longrightarrow}{_{_{n_1,n_2\rightarrow \infty}}} 0 .
%\end{equation}
%Therefore for every fixed $\lambda_0 \in S^1$ the sequence
%$a_n (\lambda_0)$ tends to a certain element $a(\lambda_0 )\in \B$.
Let $f$ be the function given by
$$
f(\lambda )= \langle \lambda  a\xi ,\eta \rangle .
$$
Then we have
$$
\abs{f_{n} (\lambda )- f(\lambda )} =
\abs{\langle (\lambda a_{n} - \lambda a )\xi ,\eta \rangle}\le
\Vert a_{n} - a\Vert \stackrel{\longrightarrow}{_{_{n\rightarrow \infty}}}0
$$
which means that the sequence $f_n$ of continuous functions
tends uniformly  to $f$. Thus $f$ is continuous and
therefore $f\in L^2_\mu (\hat{\Gamma})$ ($\mu$ is the Haar measure of the compact group $\hat{\Gamma}$).
\\
Let
$$
f (\lambda ) = \sum_{g\in \Gamma} \gamma_g \lambda_g
$$
where the righthand part is the Fourier series of $f$. Since $f_n \to f$
(in $L^2_\mu (\hat{\Gamma})$) it follows that
\begin{equation}
\label{>10}
\gamma^{(n)}_g \to \gamma_g\ \  \ {\rm  for\ \ every }\ \     g\in \Gamma
\end{equation}
where $\gamma^{(n)}_k$ are those defined by \eqref{flambda}.
Now note that property $(*)$ implies
\begin{equation}
\label{>11}
\Vert a^{(n)}_g \Vert \to \Vert E_g(a)   \Vert \ \ {\rm  for\ \ every}\ \ g\in \Gamma.
\end{equation}
And also observe that
$$
\abs {\gamma^{(n)}_g} \le \Vert a^{(n)}_g \Vert \ \ {\rm  for\ \ every}\ \ g\in \Gamma
$$
which together with \eqref{>10}, \eqref{>11}
means that
\begin{equation}
\label{>12}
\gamma_g = 0 \ \ {\rm  for\ \ every}\ \ g \in \Gamma .
\end{equation}
Now \eqref{>12} and the continuity of $f$ implies
$$
f(\lambda )=0  \ \ {\rm  for\ \ every}\ \ \lambda \in \hat{\Gamma} .
$$
In particular
$$
f(1) =<a\xi , \eta >=0 .
$$
Thus \eqref{zero} is true, whence $a=0$.
\\
For the completion of the proof of the theorem we verify the implication $(iii)\Rightarrow (ii)$.
\\
By \ref{conditionexpect}, for the proof of $(ii)$ it is sufficient to demonstrate that
\begin{equation}\label{nienormalne}
E_x(a)u_x=0,\qquad u^*_x E_{-x}(a)=0, \qquad x\in \Gamma^+.
\end{equation}
Selecting a sequence of elements $a_n\in C_0$ such that $a_n\longrightarrow a$ and taking account of
the explicit form of $E_0(a_n^*a_n)= \sum u^*_x (a^{(n)}_x)^* a^{(n)}_xu_x+ \sum (a_{-x}^{(n)})^* a_{-x}^{(n)}$, and of the fact that,
 as $n\longrightarrow \infty$, for each fixed $x\in\Gamma^+$  we have $a^{(n)}_x\longrightarrow E_x(a)$,
 $a^{(n)}_{-x}\longrightarrow E_{-x}(a)$ and $E_0(a_n^*a_n)\longrightarrow E_0(a^*a)$, we conclude that for each
 finite subset $F\subset \Gamma^+$ we have
$$
\left\| \sum_{x\in F}u^*_x E_x(a)^* E_x(a)u_x+ \sum_{x\in F}E_{-x}(a)^* u_xu_x^* E_{-x}(a)\right\|\leq E_0(a^*a).
$$
Hence equalities \eqref{nienormalne} follow by $(iii)$. The proof is complete.
 \qed

\begin{I}

\label{regular} {\bfseries Regular representation of the crossed
product}.\ \
Now we present a
faithful representation of $\A \rtimes_\al \Gamma$ that
will be written out explicitly in terms of $\A$, $\al$,
$L$.
Keeping in mind the standard regular representations for
the known various versions of crossed products we shall
call it the {\em regular representation} of  $\A \rtimes_\al \Gamma$.
In fact the construction of this
representation has been already obtained in the proof of Theorem \ref{b..1}. Let us, however, briefly discuss it.
\\
The Hilbert space $H$ is defined to  be the  direct sum $\bigoplus_{f\in
F}H^f$ of Hilbert spaces $H^f$ where $F$ is the set of all positive linear functionals on  $\A$, and
every $H^f$ is in turn  the
 direct sum of Hilbert spaces $\bigoplus_{g\in \Gamma} H^f_g$. The spaces $H^f_g$ are generated by non-negative inner
products $\langle\,\cdot\,,\,\cdot\,\rangle_g$ on the initial algebra
$\cal A$ that are given by
$$
\langle v,u\rangle_x =f\bigl(L_x(u^*v)\bigr),
\qquad
\langle v,u\rangle_{-x}
=f\bigl(u^*\al_x(1)v\bigr), \qquad x\in \Gamma^+.
$$
The  semigroup homomorphism $U\!:\Gamma^+ \to L(H)$ is defined in such a manner that for an arbitrary fixed $x\in \Gamma^+$
and for  any finite sum $h =\bigoplus_g h_g \in H^f$,  $h_g\in H^f_g$ we have  $
U_x h =\bigoplus_g(U_xh)_g$  where
$$
(U_xh)_g =
\begin{cases}
\al_x(h_{g-x}),& {\rm if}\ \ x\leq g,\\
\al_x(1) \al_g(h_{g-x}),&  {\rm if}\ \ 0 \leq g\leq x,
\\
\al_{x-g}(1)h_{g-x},&   {\rm if}\ \ g \leq 0.
\end{cases}
$$
Finally, the  representation  $\pi\!:{\cal A}\to L(H)$ is defined in such a manner that for any $a\in\A$ the operator  $\pi(a)\!:H\to H$ leaves
invariant all  the
subspaces  $H^f_g$, $F\in F$, $g\in\Gamma$, and for $h_g\in H^f_g$  we have
$$
\pi(a)h_g =
\begin{cases}
ah_g, & \ \ g\ge 0,\\
\al_{-g}(a)h_g, &\ \ g\le 0.
\end{cases}
$$

In the process of the proof of Theorem \ref{b..1} there was
verified that the triple $(H, \pi , U)$ described above is a
covariant representation of $(\A,\Gamma^+,\al)$ (in the sense of
Definition \ref{covrep}), hence by Proposition \ref{universality}
it gives the rise to a representation $(\pi\times U)$ of $\A
\rtimes_\al \Gamma$. Moreover, the covariant representation $(H,
\pi , U)$ possesses  property (*) which can be proved by repeating
the argument from  the proof of item $(iv)$ of Theorem
\ref{crossedstructure}. Thus, by Theorem \ref{isomorphism
theorem}, $(\pi\times U)$ is a faithful representation of $\A
\rtimes_\al \Gamma$.
\end{I}

\end{document}